\newif\ifstandardtemplate\standardtemplatetrue%
\newif\ifelseviertemplate\elseviertemplatefalse%
\newif\ifspringertemplate\springertemplatefalse%
\newif\ifwileytemplate\wileytemplatefalse%
\newcommand{\mypackages}{%
  \usepackage{amssymb}
  \usepackage{amsmath}
  \usepackage{amsthm}
  \usepackage{graphicx}
  \usepackage[dvipsnames]{xcolor}
  \usepackage{xfrac}
  \usepackage{todonotes}
  \usepackage{siunitx}
  \usepackage[final]{changes}
  \definechangesauthor[name=Ignacio, color=Orange]{IRO}
  \definechangesauthor[name=David, color=blue]{DPG}
  \definechangesauthor[name=Edmund, color=red]{EBN}
  \graphicspath{{Figures/}{figures/}}
  \usepackage{enumitem}
  \setlist[enumerate,1]{label = \emph{\alph*}),
                        ref   = \theenumi.\emph{\alph*}}
}

\newcommand{\mymacros}{%
  \newcommand{\vto}{VTORCH}
  \newcommand{\take}{ACPF}
  \newcommand{\concept}[1]{\textbf{\emph{##1}}}
  \newcommand{\defined}{:=}
  \newcommand{\dev}{{\mathop{\mathrm{dev}}}}
  \renewcommand{\div}{{\mathop{\mathrm{div}}}}
  \newcommand{\mbs}[1]{\boldsymbol{##1}}
  \newcommand{\pairing}[2]{\langle{##1},{##2}\rangle}
  \newcommand{\dd}[2]{\frac{\mathrm{d} ##1}{\mathrm{d} ##2}}
  \newcommand{\pd}[2]{\frac{\partial{##1}}{\partial{##2}}}
  \newcommand{\ppdd}[3]{\frac{\partial^2{##1}}{\partial{##2}\partial{##3}}}
  \newcommand{\fd}[2]{\frac{\delta{##1}}{\delta{##2}}}
  \newcommand{\ffdd}[3]{\frac{\delta^2{##1}}{\delta{##2}\delta{##3}}}
  \newcommand{\set}[1]{\left\{##1\right\}}
  \newcommand{\trace}{{\mathop{\mathrm{tr}}}}
  \newcommand{\uptohere}{\centerline{\textcolor{blue}{\rule{6cm}{0.2cm}}}}
  \let\oldLambda=\Lambda\renewcommand{\Lambda}{\mathit{\oldLambda}}
  \let\oldGamma=\Gamma\renewcommand{\Gamma}{\mathit{\oldGamma}}
}
\newcommand{\mytitle}{A fully variational numerical method for structural topology
  optimization based on a Cahn-Hilliard model}
\newcommand{\myabstract}{\noindent We formulate a novel numerical method suitable for the solution of topology optimization problems in solid mechanics. The most salient feature of the new approach is that the space and time discrete equations of the numerical method can be obtained as the optimality conditions of a single incremental potential. The governing equations define a gradient flow of the mass in the domain that maximizes the stiffness of the proposed solid, while exactly preserving the mass of the allocated material. Moreover, we propose a change of variables in the model equations that constrains the value of the density within admissible bounds and a continuation strategy that speeds up the evolution of the flow. The proposed strategy results in a robust and efficient topology optimization method that is exactly 
mass-preserving, does not employ Lagrange multipliers, and is fully variational.}
\newcommand{\myack}{The authors would like to acknowledge the funding received from the Spanish Ministry of Science and Innovation under grant PID2021-128812OB-I00.}

\ifstandardtemplate%
\documentclass[10pt,a4paper]{article}
\mypackages%
\mymacros%
\newcommand{\mybibstyle}{unsrt}

\newtheorem{theorem}    {Theorem}[section]

\theoremstyle{definition}

\newtheorem{examplex}[theorem]{$\triangleright\;$Ejemplo}

\theoremstyle{remark}

\newtheorem{remarks}{Remarks}

\usepackage{hyperref}
\usepackage[capitalise]{cleveref}
\usepackage{multirow}
\usepackage{multicol}
\usepackage{caption}
\usepackage{subcaption}
\usepackage{booktabs}
\usepackage{placeins}

\title{\mytitle}
\author{E. Bell-Navas$^1$ \and D. Portillo$^1$ \and I. Romero$^{1,2}$}
\date{$^1$ Universidad Polit\'ecnica de Madrid,
  Jos\'{e} Guti\'{e}rrez Abascal, 2, 28006 Madrid, Spain\\[2ex]%
  $^2$ IMDEA Materials Institute, Eric Kandel,2; 28906 Madrid, Spain
  \\[2ex]\today}

\newenvironment{acknowledgements}{\section*{Acknowledgements}}{}

\begin{document}
\maketitle
\begin{abstract}
  \myabstract
  \vspace{0.5em}
\end{abstract}

\fi

\ifspringertemplate%
\documentclass[smallcondensed]{svjour3}
\mypackages%
\mymacros%

\title{\mytitle%
\thanks{the thanks}}
\author{\ldots \and Ignacio Romero \and \ldots}
\journalname{The journal name}

\institute{I. Romero \at%
  IMDEA Materials Institute, Eric Kandel 2, Tecnogetafe, Madrid 28906, Spain\\
  Universidad Polit\'ecnica de Madrid, Jos\'e Guti\'errez Abascal, 2, Madrid 29006, Spain\\
  \email{ignacio.romero@imdea.org}
  \and
  XX \at%
  XXX\\
  \email{XXX}}

\titlerunning{\mytitile}
\authorrunning{I. Romero}

\date{Received: date / Accepted: date}

\maketitle
\begin{abstract}
  \myabstract%
\end{abstract}
\fi

\ifelseviertemplate%
\documentclass[preprint,11pt]{elsarticle}
\mypackages%
\mymacros%
\newcommand{\mybibstyle}{unsrt}

\journal{Computer Methods in Applied Mechanics and Engineering}

\begin{frontmatter}

\title{\mytitle}

\author[1,2]{Ignacio Romero\corref{cor1}}
\ead{ignacio.romero@upm.es}

\author[3]{xxx}
\ead{xxx}

\address[1]{Universidad Polit\'ecnica de Madrid, Spain}
\address[2]{IMDEA Materials Institute, Spain}
\address[3]{xx}

\cortext[cor1]{Corresponding author. ETS Ingenieros Industriales,
  Jos\'e Guti\'errez Abascal, 2, Madrid 28006, Spain}

\begin{abstract}
\myabstract%
\end{abstract}

\end{frontmatter}

\fi

\ifwileytemplate%
\documentclass[doublespace,times]{nmeauth}
\mypackages%
\mymacros%
\usepackage{natbib}
\message{Compiling with Wiley template}

\runningheads{I. Romero}{\dots}

\title{\mytitle}
\author{Ignacio Romero\affil{1}\affil{2}\corrauth}

\address{\affilnum{1} ETSII, Universidad Politécnica de Madrid,
         Jos\'{e} Guti\'{e}rrez Abascal, 2, 28006 Madrid, Spain\break%
         \affilnum{2} IMDEA Materials Institute, Eric Kandel 2, 28096 Getafe, Madrid, Spain}

\corraddr{Dpto.~de Ingenier\'{\i}a Mec\'{a}nica; E.T.S. Ingenieros Industriales;
José Gutiérrez Abascal, 2; 28006 Madrid; Spain. Fax (+34) 91 336 3004}

\begin{abstract}
  \myabstract
\end{abstract}

\maketitle
\fi


\section{Introduction}
\label{sec-intro}
In the context of solid mechanics, basic topology optimization (TO) algorithms strive to find the lightest
structure or solid that can resist some given loading without exceeding certain deformation. Equivalently, these methods try to find the structure of minimum compliance for a given fixed volume. This is a fundamental problem in mechanical design that has motivated the development of a large number of alternative methodologies in an attempt to
solve it (see, e.g., the monographs \cite{allaire2002to,bendsoe2004ik} or the reviews
\cite{deaton2013uj,yago2021cd,sigmund2013ht,mukherjee2021hq} and the references therein).

The most widely employed TO methods are written as optimization programs in which the cost functional is a certain measure of the compliance or the work done by the external forces on the structure, and the constraints include the total volume of the solid. The optimality conditions of this problem lead to the computation of the sensitivities
and, ultimately, to iterative solutions in which the solid's mass gets transported to the regions where it can contribute most to the structural stiffness. It is well-known that these programs are ill-posed \cite{sigmund1998ill,lazarov2016length}, and numerical methods devised to approximate their solutions must, by force, resort to regularization or
stabilization. There is no canonical way to carry out this regularization and this explains the emergence of quite 
different families of methods such as those based on homogenization~\cite{bendsoe1988ut}, phase-field approaches based on the solution of the Allen-Cahn problem
(\cite{takezawa2010ps,gao2020oc,bartels2016ac,marino2021}, thermodynamic-like formulations
\cite{junker2015pn,junker2021iv}, perimeter constraints~\cite{tavakoli2014ie}, variational methods for the Cahn-Hilliard 
equation~\cite{bartels2021ch}, etc. This variety of solution strategies points to the fact that there is still room for new approximation schemes that combine robustness and efficiency.

In the current work, we follow the last of the techniques identified above, transforming the TO problem into an ancillary Cahn-Hilliard dynamical system. In these types of formulations, the main unknown is the volumetric fraction, also referred to as the ``relative density'' or just the ``density'', a scalar field that indicates the material distribution within the body. In particular, Cahn-Hilliard formulations cast the topology optimization problem into a gradient flow \cite{blank2012,kim2016}, with fluxes redistributing existing mass while maximizing, for example, the stiffness of the solid. The numerical solution of this evolution problem is not exempt from difficulties, especially due to the constraint that the density values must remain between~0 and~1.

One particularly powerful tool for designing space and time approximations for a general class of initial boundary value problems is provided by the so-called \emph{variational updates} in combination with Galerkin projections. These methods propose an incremental functional whose stationarity conditions are consistent approximations of the time-discretized evolution equations of the problem. Initially proposed for purely mechanical problems \cite{ortiz1999tq,radovitzky1999kc}, they have
proven quite successful in the solution of complex thermo-mechanical problems
\cite{yang2006vg,canadija2011hr}, microstructural evolution \cite{carstensen2002uh},
thermo-chemo-mechanical problems \cite{romero2021dd}, hydrogen diffusion \cite{andres2024on}, and
topology optimization \cite{bartels2016ac,bartels2021ch}.

Variational updates possess several interesting properties. First, in some situations, the variational nature of the problem statement allows for deep mathematical analyses of existence based on Gamma-convergence \cite{schmidt2008yv}. From the numerical standpoint, these are implicit methods that have always
symmetric tangents, in contrast with standard discretizations for multi-field problems which
typically yield non-symmetric tangents (see \cite{romero2021dd} for the discussion). Simply because of this property, the memory requirements for variational methods are (asymptotically) 50\%\ smaller than
for standard implementations, and solution times for direct solvers are also around half of the standard counterparts.

While the reasons adduced above should suffice to motivate the formulation of a variational method for topology optimization, we must point out that any deviation from the variational program would spoil the symmetry of the tangent operator. For example, whereas the thermodynamic approach of Junker and coauthors \cite{junker2015pn,junker2021iv} has a variational foundation, the time discretization breaks the optimality of the formulation, meaning that the final equations are not the stationarity conditions of a single functional and, hence, the symmetry of the tangent is lost. Similarly, a variational method was recently proposed \cite{bartels2021ch} which is closely related to the present article. However, in this work, the constraint that the values of the density field must remain within 0 and~1 is
imposed in a non-variational fashion, and the symmetry of the tangent operator is spoiled.

The work described next is, to the authors' knowledge, the first numerical method for topology optimization that is fully variational. The derivations that follow will serve to identify a single
incremental functional whose stationarity conditions, as advanced, will coincide with the update equations for all the fields in the solution. Naturally, the tangent operator will retain its symmetry in the discrete setting and this can be exploited by special, fast linear solvers. One additional novelty of our proposed formulation is the elimination of the constraint that the density values should remain in $[0,1]$. Instead, we will advocate for a change of variable such that the density will always remain in this interval. In our formulation, the unknown field will be unconstrained and mapped onto $(0,1)$ to recover the sought density field. To speed up computations further, a continuation method will be proposed that modifies the Cahn-Hilliard dynamics to favor the system's evolution towards a state of maximum potential energy that will coincide with the topologically optimal solution.

The remainder of the article has the following structure. In Section~\ref{sec-statement}, we
summarize the problem of topology optimization for small-strain solid mechanics and
reformulate it as a Cahn-Hilliard gradient flow. Building on this formulation, we
identify in Section \ref{sec-var-update} an incremental variational problem whose stationarity conditions are consistent
approximations of the evolution equations previously identified. Numerical simulations of topology optimization problems are shown in
Section~\ref{sec-examples}, where the new method is compared with existing ones.
Section~\ref{sec-summary} closes the article with a summary of the main results.

Notation: Throughout the article, vector- and tensor-valued variables and functions will be indicated in boldface. Points in $\mathbb{R}^d$ with $d\ge2$ will also be written in boldface.

\section{Problem statement}
\label{sec-statement}
We first present a model for topology optimization. As we will later explain, the basic statement of the problem is the starting point for many different models that perturb it to make it more tractable. In this section, we will proceed step by step to justify the one we will favor. We advance that most decisions and simplifications are made to set the generic topology optimization problem in a ``nice'' functional setting, one that is amenable to numerical discretization.

\subsection{Shape optimization}

We consider a topology optimization problem that aims to maximize the potential energy of an elastic body in equilibrium under known boundary conditions and loading, and under the constraint that the total measure of the domain must remain below a given bound. To define this problem precisely, let us consider first the classical problem of elasticity for a fixed domain. Let $\Omega$ be a connected open set in $\mathbb{R}^d$, with $d=2$ or $3$, and a smooth boundary $\partial\Omega$ that can be split into three disjoint sets, namely $\partial_D\Omega$, $\partial_0\Omega$ and $\partial_t\Omega$. We are interested in finding the displacement field $\mbs{u}$ in the space
\begin{equation}
  \label{eq-space-S}
  \mathcal{S}(\Omega) :=
  \left\{
  \mbs{v} \in [H^1(\Omega)]^d\ ,\; \mbs{v}(\mbs{x}) = \mbs{0}\ \text{for}\ \mbs{x}\in\partial_D\Omega
  \right\}
\end{equation}
that minimizes the potential energy
\begin{equation}
  \label{eq-potential-energy}
  V(\mbs{v};\Omega)
  \defined
  \int_{\Omega} W(\mbs{\varepsilon}(\mbs{v}); \mbs{x}) \; \mathrm{d} V -
  \int_{\partial_t\Omega} \mbs{t}(\mbs{x})\cdot \mbs{v}(\mbs{x}) \; \mathrm{d} A
\end{equation}
among all the functions in $\mathcal{S}(\Omega)$. In these equations, $[H^1(\Omega)]^d$ refers to the Hilbert space of vectors fields defined on $\Omega$ whose components are square-integrable and also their first derivatives; the stored energy function $W$ depends on the strain
$\mbs{\varepsilon}=(\nabla \mbs{u}+\nabla^T\mbs{u})/2$ and the function itself can change from point to
point. The tractions $\mbs{t}$ are assumed to be non-zero almost everywhere on $\partial_t\Omega$. The set $\partial_0\Omega$ is also part of the Neumann boundary but is free of tractions and plays a different role in the topology optimization problem; thus, we consider it separately from $\partial_t\Omega$.

The topology optimization that we are interested in considers all the bodies whose domains belong to the set
\begin{equation}
  \label{eq-b-alpha}
  \begin{split}
  \mathcal{B}_{\alpha}
  \defined
  \big\{
  & \omega\subseteq\Omega\ , \ \partial\omega = \partial_D\omega\cup\partial_t\omega\cup\partial_0\omega\ , \\
  & \partial_D\omega = \partial_D\Omega, \
    \partial_t\omega = \partial_t\Omega,\
    |\omega| = \alpha\,|\Omega|\ ,
  \big\}\
  \end{split}
\end{equation}
where $|\omega|$ denotes the measure of any set $\omega\subseteq\Omega$. More precisely,
given some $\alpha\in(0,1]$, the topology optimization problem seeks the optimal domain $\omega$ that solves
\begin{equation}
  \label{eq-topopt}
  \sup_{\omega\in \mathcal{B}_{\alpha}}\inf_{\mbs{u}\in \mathcal{S}(\omega)} V(\mbs{u};\omega)\ .
\end{equation}
Strictly speaking, elements of $\mathcal{B}_{\alpha}$ cannot be arbitrary subsets of $\Omega$ and must possess
some kind of regularity so that problem~\eqref{eq-topopt} is well-posed. The
problem statement~\eqref{eq-topopt} is very general and difficult to study, both analytically as well as numerically. Thus, we will refrain from attempting to solve it directly, and important questions such as the topology of $\mathcal{B}_{\alpha}$ will not be addressed. We refer to the monograph by Bucur and Butazzo \cite{bucur2005me} for a deep analysis of problems of this type.

\subsection{Reparameterized formulation}
\label{subs-relaxed}
Problem~\eqref{eq-topopt} might not have a solution. If the material properties satisfy standard bounds, the infimum of $V$ for a fixed domain $\omega$ is actually a unique minimum that
always exists, provided that $\omega$ satisfies some regularity conditions. However, the dependency of the functional on the set $\omega$ is complex, and the space $\mathcal{B}_{\alpha}$ itself is not even linear. Given these difficulties, we choose to
replace problem~\eqref{eq-topopt} with a simpler one, one that might enable some kind of numerical approximation.

This strategy demands further justification. We depart from problem~\eqref{eq-topopt} due to its complexity. We postulate a new problem that incorporates our goal for topology optimization but do not claim that the solution to the latter converges, in any sense, to the solution of the former, which might not even exist. Restricting the problem statement might inadvertently lose an optimal solution, but we accept it in search of a more flexible formulation, one that is better suited to approximation.

Thus, instead of searching for an optimal domain in the set $\mathcal{B}_{\alpha}$, we settle for finding a \emph{density} field $\rho$ with values in $[0,1]$ that can be interpreted as the ratio of material at a given point of $\Omega$ (values $\rho(\mbs{x})=0$ indicate that there is no material at point $\mbs{x}$, whereas
$\rho(\mbs{x})=1$ means that the point $\mbs{x}$ is a point of the solution shape). This \emph{phase-field} approach allows for intermediate values $\rho(\mbs{x})\in(0,1)$ that, while nonphysical, open
the door to optimizing some functional among the objects of a space with ``nice'' properties. 

The functional~\eqref{eq-potential-energy} needs to be adapted to the new setting in a way that solids with optimal mass distributions are found. For example, and slightly abusing the notation,
problem~\eqref{eq-topopt} could be written as
\begin{equation}
\label{eq-opt-rho}
\sup_{\rho}\inf_{\mbs{u}\in \mathcal{S}(\Omega)}
V(\mbs{u};\rho)\ ,
\end{equation}
with
\label{eq-opt-rho2}
\begin{equation}
V(\mbs{u};\rho)
:=
  \int_{\Omega} W(\mbs{\varepsilon}(\mbs{v}); \mbs{x})\;\rho(\mbs{x}) \; \mathrm{d} V -
  \int_{\partial_t\Omega} \mbs{t}(\mbs{x})\cdot \mbs{v}(\mbs{x}) \; \mathrm{d} A
\ ,
\end{equation}
where we have intentionally left undefined the space to which the density field belongs.

As we will see, the most convenient functional setting for a topology optimization model is that of
Hilbert spaces, since it will naturally suggest a spatial discretization based on the
Galerkin method. In the case of the density $\rho$, it must also 
belong to $L^\infty(\Omega)$ so that we can constraint its values (almost everywhere) 
to the interval $[0,1]$ by imposing $\|2\rho-1\|_{L^\infty(\Omega)}\le 1$. 

To sidestep the difficulties involved with the constraints on the density values, we choose to reformulate further the TO problem in terms of a pseudo-density $\theta$ defined on the space
\begin{equation}
  \mathcal{T}
  \defined
\left\{
\theta \in H^1(\Omega),\; \theta(\mbs{x}) = \bar{\theta}(\mbs{x})\ \text{for}\ \mbs{x}\in\partial_D\Omega\cup\partial_t\Omega
\right\}\,,
\end{equation}
for some known $\bar{\theta}$, and then recover the density field by selecting a suitable smooth, bijective map $\hat{\rho}:\mathbb{R}\to(0,1)$ such that
$\rho=\hat{\rho}(\theta)$. For example, arbitrary density fields can be obtained by composing functions in $\mathcal{T}$ with a scaled \emph{logistic} function $L_k:\mathbb{R}\to(0,1)$ defined as
\begin{equation}
\label{eq-logistic}
L_k(z) \defined \frac{1}{1+e^{-k\,z}} ,
\end{equation}
where $k>0$ is a scale parameter. By construction, the mapping $\hat{\rho}\equiv L_k$ satisfies all the desired properties and can be used to recover the density from the pseudo-density.

To accommodate the constraint that the total volume of the optimal solution must be identical to the available volume, a constraint  intrinsic to the definition of the set $\mathcal{B}_{\alpha}$ (see Eq.~\eqref{eq-b-alpha}), we introduce the ancillary sets
\begin{equation}
\mathcal{T}_\alpha
\defined
\{ \theta\in \mathcal{T},\;
\int_\Omega \hat{\rho}(\theta(\mbs{x}))\,\mathrm{d}V = \alpha\,|\Omega|
\} ,
\end{equation}
that constrains further the admissible set for the quasi-densities~$\theta$.

With the previous definitions, we can already formulate a TO problem that will attempt to obtain
densities with values in $(0,1)$ and where the total mass of the solution is constrained to a given value. However, if the proposed solution is to be similar to the one provided by the original TO problem, the density
field should resemble the characteristic function of an optimal set. For that, intermediate values of the
density, i.e., far from 0 and 1 should be penalized in the functional. To promote that the density clusters around the values $\{0,1\}$, while avoiding the appearance of fine structures, we introduce the Modica-Mortola functional \cite{marino2021}:

\begin{equation}
\label{eq-modica}
M_\epsilon(\theta)
\defined
\gamma \int_\Omega
\left[
\frac{U(\theta(\mbs{x}))}{\epsilon} + \frac{\epsilon}{2} |\nabla\theta(\mbs{x})|^2
\right]\,\mathrm{d}V,
\end{equation}
with $\gamma$ being a constant with dimensions of energy over area, $\epsilon>0$ a regularization
parameter, and $U:\mathbb{R}\to\mathbb{R}_{\ge0}$ defined to be the double well $U(z) \defined
2 z^2(1-z)^2$. Then, the desired goal of keeping the density close to 0 and 1 can be attained simply by appending $M_{\epsilon}$ to the definition of the functional $V$ and taking $\epsilon\ll L$,
where $L$ is a characteristic length of the domain $\Omega$; its diameter, for example.

Collecting all the previous definitions, we can postulate a TO problem that will find the optimal mass distribution $\rho$ and its corresponding displacement field as the solutions to
\begin{equation}
  \label{eq-prob-ve}
  \rho = \hat{\rho}(\theta)\,,
  \qquad
  (\theta,\mbs{u}) =
  \arg
  \sup_{\theta\in \mathcal{T}_{\alpha}} \inf_{\mbs{u}\in \mathcal{S}(\Omega)}
  F_{\epsilon}(\mbs{u},\theta)\ ,
\end{equation}
where we have introduced the regularized energy
\begin{equation}
  F_{\epsilon}(\mbs{u},\theta)
  \defined
  V(\mbs{u},\hat{\rho}(\theta))  -   M_\epsilon(\theta)\ .
\end{equation}

It is much easier to work with program~\eqref{eq-prob-ve} than with~\eqref{eq-topopt}. Still, the definition of the space $\mathcal{T}_{\alpha}$ includes a \emph{global} constraint on the total mass of the body.
These types of constraints are inconvenient from the numerical standpoint since they couple all the degrees of freedom of a discretized model. In anticipation of the methods that we will later introduce, in Section~\ref{subs-evolution} we will propose a convenient way of removing the volume constraint from the problem.

\begin{remarks}\
\par
\begin{enumerate}

\item The density fields that result from the proposed optimization can never attain the values $\{0,1\}$, but can be as close to these values as necessary. This is a limitation of the change of variable from $\theta$ to $\rho$, and it cannot be avoided if we insist that this map should be a smooth bijection.

\item The choice for the logistic function again is not unique. It can be replaced with any smooth and bijective mapping from the real line to $(0,1)$. In our choice, there is a shape parameter $k$ that will be adjusted for convenience but could be fixed to an arbitrary non-zero number.

\item The expression selected for the double well is somewhat arbitrary too. The function $U$ should
  be differentiable and have two minima when $\rho\approx0$ and $\rho\approx1$.

\end{enumerate}
\end{remarks}

\subsection{Topology optimization as an evolution problem}
\label{subs-evolution}
To sidestep the global constraint that appears in the definition of the TO
program~\eqref{eq-prob-ve}, we propose to shift the structure of the problem from a saddle-point
formulation to a dynamical system, one designed to possess an attractor that coincides with the sought optimum, and whose flow preserves the total mass. For that, and abusing the notation, let us now look for a density $\rho\in C^1([0,T],\mathcal{S}(\Omega))$ that evolves driven by a time-dependent vector field of \emph{mass flux} denoted as $\mbs{j}$, as given by the transport problem
\begin{subequations}
  \label{eq-transport}
    \begin{align}
      \pd{\rho}{t}(\mbs{x},t) + \nabla\cdot \mbs{j}(\mbs{x},t)
      &= 0\ ,&\qquad \mbs{x}\in \Omega\ , \\
      \mbs{j}(\mbs{x},t)\cdot \mbs{n}(\mbs{x}) &= 0\ , &\qquad \mbs{x}\in \partial\Omega\ ,
      \label{eq-transport2} \\
      \rho(\mbs{x},0) &= \rho_0(\mbs{x})\ , &\qquad \mbs{x}\in \Omega\ ,
    \end{align}
\end{subequations}
with $\rho_0$ being an initial density field satisfying
\begin{equation}
    \label{eq-initial-rho}
    \int_{\Omega}\rho_0(\mbs{x}) \; \mathrm{d} V = \alpha\,|\Omega|\ .
\end{equation}
In Eq.~\eqref{eq-transport}, the field $\mbs{n}$ refers to the outward, unit vector normal to the boundary $\partial\Omega$, and $\nabla\cdot$ is
the divergence operator. Eqs.~\eqref{eq-transport} describe the evolution of a density field in
$\Omega$ that, since there is no mass influx through the boundary nor through any source, must satisfy
\begin{equation}
  \label{eq-proof-cm}
  \dd{}{t} \int_{\Omega} \rho(\mbs{x},t)\; \mathrm{d} V =
  - \int_{\Omega} \nabla\cdot \mbs{j}(\mbs{x},t) \; \mathrm{d} V = 0 .
\end{equation}
Let us postulate an expression for the mass flux such that
when $t\to\infty$ the density field will define a body with maximum elastic
energy. For that, let us introduce first the \emph{pseudo chemical potential} $\mu$ defined as
\begin{equation}
  \label{eq-chemical-pot}
  \mu(\mbs{x}) \defined
  \frac{1}{\hat{\rho}'(\theta(\mbs{x}))}
  \fd{G_\epsilon}{\theta}(\theta(\mbs{x}))\ ,
  \qquad\text{with}\qquad
  G_{\epsilon}(\theta) \defined \inf_{\mbs{u}\in \mathcal{S}(\Omega)} F_{\epsilon}(\mbs{u},\theta)\ .
\end{equation}
Also, assume we select a \emph{convex} kinetic potential $\psi: \mathbb{R}^d\to \mathbb{R}$, a function of $\nabla\mu$, such that
\begin{equation}
\label{eq-kinetic}
\psi(\mbs{0})=0,
\qquad
\pd{\psi}{\nabla\mu}(\mbs{0})=\mbs{0}, \qquad
\mbs{j} = \pd{\psi}{\nabla\mu}(\nabla\mu)\ .
\end{equation}
Then,
\begin{equation}
  \label{eq-proof}
  \begin{split}
    \dd{}{t} G_\epsilon(\theta)
    &=
      \int_{\Omega} \fd{G_\epsilon}{\theta}(\theta)\; \dot{\theta}(\mbs{x}) \; \mathrm{d} V\\
    &=
      \int_{\Omega}
      \frac{1}{\hat{\rho}'(\theta(\mbs{x}))}
      \fd{G_\epsilon}{\theta}(\theta)\; \dot{\rho}(\mbs{x}) \; \mathrm{d} V\\
    &=
       \int_{\Omega} \mu(\mbs{x})\; (-\nabla\cdot \mbs{j}(\mbs{x})) \; \mathrm{d} V \\
    &=
       \int_{\Omega} \nabla\mu(\mbs{x}) \cdot \partial\psi(\nabla\mu(\mbs{x})) \; \mathrm{d} V
      -
      \int_{\partial_D\Omega} \mu(\mbs{x})\; \mbs{j}(\mbs{x})\cdot \mbs{n}(\mbs{x}) \; \mathrm{d} A
    \\
    &\geq 0\ ,
  \end{split}
\end{equation}
where we have employed the convexity of $\psi$, relations~\eqref{eq-kinetic}, and the boundary
condition~\eqref{eq-transport2}. For simplicity, in the simulations of Section~\ref{sec-examples},
the kinetic potential will be taken as a quadratic function of the norm of $\nabla\mu$. The dynamical properties that we seek are a consequence of the convexity of this potential and any other that shares this property will generate
a flow that will get closer and closer to a maximizer of~$G_{\epsilon}$.

To sum up, the topology optimization problem can be thought of as a mass transport problem in which the
mass flux is driven by the maximization of $G_\epsilon$, while the displacement field $\mbs{u}$ is
the one that ensures mechanical equilibrium at every time instant.
In equilibrium thermodynamics, the chemical potential can be obtained as the partial derivative of
the grand canonical potential with respect to the density. In the current model, the role of the grand canonical potential is played by~$G_\epsilon$ (see Eq.~\eqref{eq-chemical-pot}). Strictly speaking, however,
thermodynamics is never invoked to obtain $G_\epsilon$. Instead, we view this function simply as a flow generator that possesses the right properties.

\section{A variational method for the evolution problem of
topology optimization}
\label{sec-var-update}
We propose now a numerical method based on the concept of \emph{variational updates}
\cite{ortiz1999tq,romero2021dd} that identifies a discrete \emph{incremental potential} whose variations
with respect to the unknown fields are first-order approximations (in time) to the evolution
equations of TO obtained in Section~\ref{sec-statement}.

\subsection{Time semidiscretization}
To start the time discretization, let us first consider a partition of the time interval of interest into subintervals $(t_k,t_{k+1})$ with
$k=0,1,\ldots,M$ and $0\equiv t_0<t_1<\cdots<t_M \equiv T$ and $\Delta t_k\defined t_{k+1}-t_k$.
Then, the approximation to the fields $\mbs{u},\theta,\mu$ at time $t_{k+1}$ will be denoted, respectively, as $\mbs{u}_{k+1},\theta_{k+1},\mu_{k+1}$ and obtained from the stationarity conditions of the potential
\begin{equation}
  \label{eq-incremental}
  \begin{split}
    J_k(\mbs{u}_{k+1},\theta_{k+1},\mu_{k+1}) \defined&
                                           F_{\epsilon}(\mbs{u}_{k+1},\theta_{k+1}) - F_{\epsilon}(\mbs{u}_k,\theta_k) \\
    &- \int_{\Omega}   \left[
      (\hat{\rho}(\theta_{k+1})-\hat{\rho}(\theta_k))\;\mu_{k+1} + \Delta t_k\, \psi(\nabla \mu_{k+1}) \right] \; \mathrm{d} V\ ,
  \end{split}
\end{equation}
where $(\mbs{u}_k,\theta_k,\mu_k)$ are assumed to be known fields. Indeed, the
stationarity conditions of this functional are
\begin{subequations}
  \label{eq-stationarity}
  \begin{align}
    0 &= D_1 J_k\cdot \mbs{v}
        = \int_{\Omega}
        \pd{W}{\mbs{\varepsilon}}(\mbs{\varepsilon}(\mbs{u}_{k+1}(\mbs{x})))\cdot
        \nabla^s \mbs{v}(\mbs{x})\;
        \hat{\rho}(\theta_{k+1}(\mbs{x})) \; \mathrm{d} V -
        \int_{\partial_t\Omega}
        \mbs{t}(\mbs{x})\cdot \mbs{v}(\mbs{x}) \; \mathrm{d} A \ ,
        \label{eq-stationarity-1}\\
    0 &= D_2 J_k\cdot \zeta =
        \int_{\Omega} \left[
        \fd{F_{\epsilon}}{\theta}(\mbs{u}_{k+1}(\mbs{x}),\theta_{k+1}(\mbs{x}))
        - \hat{\rho}'(\theta_{k+1}(\mbs{x})) \mu_{k+1}(\mbs{x})
        \right] \zeta(\mbs{x}) \; \mathrm{d} V\ ,
        \label{eq-stationarity-2}\\
    0 &= D_3 J_k\cdot \eta =
        - \int_{\Omega} \left[
        (\hat{\rho}(\theta_{k+1})-\hat{\rho}(\theta_k))\,\eta(\mbs{x}) +
        \Delta t_k\, \partial\psi(\nabla\mu_{k+1}(\mbs{x})) \cdot \nabla\eta(\mbs{x})
        \right]  \; \mathrm{d} V\ ,
        \label{eq-stationarity-3}
  \end{align}
\end{subequations}
where $\mbs{v}\in \mathcal{S}(\Omega), \zeta\in H^1(\Omega), \eta\in H^1(\Omega)$ are arbitrary functions.

Eq.~\eqref{eq-stationarity-1} expresses the mechanical equilibrium in the body with density field
$\hat{\rho}(\theta)$, when both the displacement and the pseudo-density are evaluated at $t_{k+1}$.
Eq.~\eqref{eq-stationarity-2} coincides with the definition of chemical potential and recovers
Eq.~\eqref{eq-chemical-pot} when all the fields are evaluated at $t_{k+1}$. Last, Eq.~\eqref{eq-stationarity-3} is the weak statement of the mass transport equation~\eqref{eq-transport}, discretized in time with a first-order scheme similar to the backward-Euler implicit method. 
The three algebraic equations in~\eqref{eq-stationarity} correspond to the time discretization of a differential-algebraic system in which only the transport equation is differential. Notice that the rate approximation is not exactly the difference of the unknown field at consecutive time instants,
hence the proposed time integrator is not identical to the implicit Euler method. This seemingly innocuous
difference, which does not add any overhead to the computation of the discrete evolution equations, renders the second linearization symmetric, guaranteeing all the computational advantages described in Section~\ref{sec-intro}.

\subsection{A continuation method for the double-well definition}
\label{sub-cont-doublewell}

The double-well potential in the Modica-Mortola functional causes numerical difficulties for the
convergence of iterative solvers for TO when starting from a uniform density distribution.
In these situations, very small time steps are required to trigger phase
separation. To address this problem and allow for larger time step sizes during the initial phases of the solution, a continuation method is presented below.

The proposed method starts the simulation using a degenerate double-well potential that simplifies the convergence of the solver when the density is homogeneous. The basins of the double-well are progressively split apart, ensuring that the phase separation is smooth and that a complete separation is guaranteed at the end of the simulation. In this way, the double-well becomes time-dependent, i.e., of the form
\begin{equation}
\label{eq-doubleWellTime}
U(\theta,t) = 2 (\theta-\theta^{-}(t))^2 (\theta-\theta^{+}(t))^2\ .
\end{equation}
Using this expression for the potential, or a similar one, its minima depend on time, effectively accomplishing the objective of 
separating its wells as the simulation progresses.

We propose here a change of variable based on the logistic function~\eqref{eq-logistic}, employing
$\theta^{{\pm}}(t) = \pm \Delta \theta(t)$, with $\Delta \theta$
being the distance from the origin (where the local maximum of the double-well is located) to the
minima. The time evolution of $\theta^{\pm}$ should smooth the transition from the uniform initial condition to the final distribution. A linear expression is proposed here for the change of variable used up to a maximum that is obtained in a characteristic time $T$,
\begin{equation}
\label{eq-continuationTheta}
\Delta \theta(t) = \min \left( (\Delta \theta_T - \Delta \theta_0)\frac{t}{T} ,\; \Delta \theta_T \right)
\end{equation}
where $\Delta \theta_0$ and $\Delta \theta_T$ are the values of $\Delta \theta$ at the beginning of the simulation and at time $T$ respectively. The value of $\Delta \theta_0$ should be selected such that convergence at the beginning of the simulation is easily attained, e.g., $\Delta \theta_0 = \theta_0$. Likewise, $\Delta \theta_T$ must be chosen such that, at the end of the simulation, the double-well minima correspond to density values close to $0$ and $1$ (typically $\rho_{\min} = 10^{-3}$).

\begin{figure}[t]
    \centering
    \includegraphics[width=0.8\textwidth]{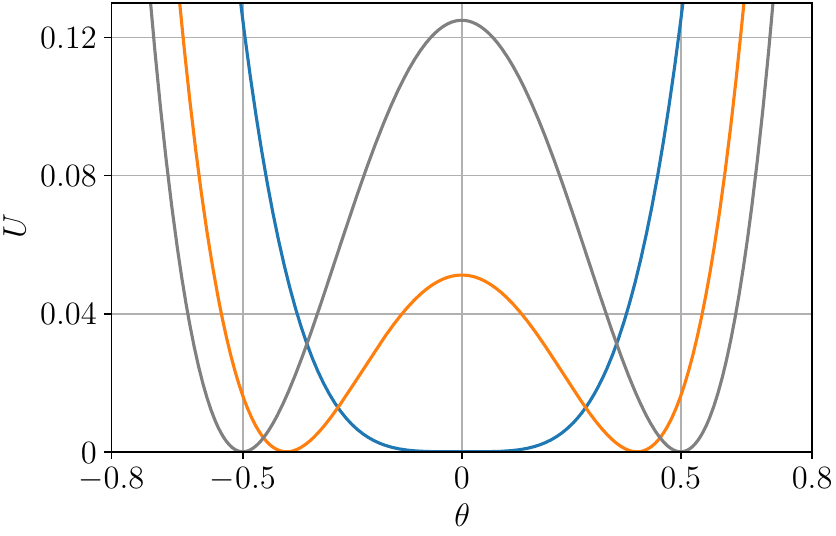}
    \caption{Continuation method for the double-well potential. Potential at times $t/T=0$ (blue), $t/T=0.8$ (orange), and $t/T=1$ (gray).}
    \label{fig-continuation}
\end{figure}

Figure \ref{fig-continuation} illustrates the behavior of the double well function $U$ as
$\theta^{\pm}(t)$ changes over time in a problem in which the initial condition is $\rho_0 = 0.5$
($\theta_0 = 0$) and $\Delta\theta(t)=0.5 \cdot t$. Initially, the double-well potential degenerates to a function with a
single minimum; as the simulation progresses, the wells slowly split apart, finally developing the desired
minima at $\theta=\pm 0.5$.

Using the continuation method proposed, the evolution of the mass redistribution process exhibits
two different phases:
\begin{itemize}
\item The first steps of the problem are characterized by a fast transport of mass. In this stage,
  the double-well potential is relatively ``flat'' and density is allowed to take values far from 0 and~1, while mass
  rapidly flows within the solution domain.

    \item In later times, the imposed double-well potential develops two minima, density values far from 0 and 1 are
      penalized, and a structure with two distinctive phases emerges. In this phase, mass
      flows only to sharpen the interfaces, eventually reaching a stationary solution.
\end{itemize}

\subsection{Spatial discretization}
Eqs.~\eqref{eq-stationarity} describe a time semidiscretization of a model for topology
optimization. A complete spatial \emph{and} temporal discretization may be accomplished by
restricting the unknown fields to finite element spaces, for example. To introduce the latter, consider a triangulation of $\Omega$ into disjoint elements $\{\Omega_e\}_{e\in \mathcal{E}}$ where $\mathcal{E}$ is the \emph{element set}. These elements are assumed to be a regular partition of $\Omega$ and are connected at nodes with labels $a\in \mathcal{N}$. To define the finite element approximation of Eq.~\eqref{eq-incremental} we need to introduce finite-dimensional subspaces of $\mathcal{S}(\Omega)$ and $H^1(\Omega)$ which we denote, respectively, as $\mathcal{S}_h$ and $\mathcal{T}_h$. These two are given by
\begin{equation}
  \label{eq-fem-spaces}
  \begin{split}
    \mathcal{S}_h &=
                    \left\{
                    \mbs{v}^h\in \mathcal{S}(\Omega), \
                    \mbs{v}^h(\mbs{x}) = \sum_{a\in \mathcal{N}} N^a(\mbs{x})\;\mbs{v}_a,
                    \ \mbs{v}_a\in \mathbb{R}^d
                    \right\} ,\\
    \mathcal{T}_h &=
                    \left\{
                    {\zeta}^h\in H^1(\Omega), \
                    {\zeta}^h(\mbs{x}) = \sum_{a\in \mathcal{N}} N^a(\mbs{x})\;\zeta_a,
                    \ \zeta_a\in \mathbb{R}
                    \right\} ,\\
  \end{split}
\end{equation}
where $N^a:\Omega\to \mathbb{R}$ are standard finite element shape functions. By restricting the functional~\eqref{eq-incremental} to $\mathcal{S}_h\times \mathcal{T}_h \times \mathcal{T}_h$ we find that the final approximation problem consists of
solving
\begin{equation}
  \label{eq-space-time-discrete}
  (\mbs{u}_{k+1}^h,\theta_{k+1}^h,\mu_{k+1}^h)
  =
  \arg \max_{\zeta^h\in\mathcal{R}_h}\min_{\mbs{v}^{h}\in \mathcal{S}_h, \eta^h\in \mathcal{R}_h}
  J_k(\mbs{v}^h, \zeta^h, \eta^h)\ ,
\end{equation}
when $(\mbs{u}^h_k, \theta^h_k, \mu^h_k)\in \mathcal{S}_h\times \mathcal{R}_h\times \mathcal{R}_h$ are
given. The solutions to the optimization program \eqref{eq-space-time-discrete} provide a discrete
semigroup $G_{\Delta t}$ defined by $(\mbs{u}^h_k, \theta^h_k, \mu^h_k)\mapsto(\mbs{u}^h_{k+1},
\theta^h_{k+1}, \mu^h_{k+1}) = G_{\Delta t_k}(\mbs{u}^h_k, \theta^h_k, \mu^h_k)$ that can generate
the solution trajectory from the initial conditions $(\mbs{u}_0^h, \theta_0^h,\mu_0^h)$. At every
discrete step, the desired density field is univocally determined by the relation
$\rho_{k}=\hat{\rho}(\theta_k^h)$. We note that this density is not a member, in general, of the
finite element space but has values in~$(0,1)$ always.

In this section, we have introduced a Variational method for the Topology Optimization with
Reparameterization and Cahn-Hillard form, and henceforth, we will refer to it, for simplicity, as \vto. The algorithm
requires the implicit integration of five nonlinear, coupled equations with five degrees of freedom per node of
the finite element model.

\subsection{Parameter selection}
\label{choice-parameters}

One of the major drawbacks of topology optimization methods, regardless of their type, is that they contain parameters that usually require fine-tuning for each problem. This adjustment is often \textit{manual}, depends critically on the user's experience, and limits the robustness of the overall procedure. In particular, the method presented in this article depends on four parameters ($\epsilon$, $\gamma$, $\kappa$, $k$) that are not directly obtained from experimental measurements of the material. A rational choice of these parameters is presented in the following, based on their mathematical interpretation and dimensional analysis. Based on the foregoing remarks, the parameters can be set automatically, increasing the robustness of the approach and reducing the need for trial and error:

\begin{itemize}

    \item The parameter $\epsilon$ controls the regularization term in the Modica-Mortola
      potential~\eqref{eq-modica}. It is known to be related to the interface thickness, $\delta$, and already
      in the original article by Cahn and Hilliard~\cite{cahn1958aw}, an estimate of $\delta$ in terms of $\epsilon$ was proposed for the one-dimensional case. These estimates have been
      later refined in the literature, in particular for double-well potentials \cite{wodo2011jcp}.
      One of these relationships has the advantage of being applicable to any double-well potential,
      which is why, in this work, we employ
    \begin{equation}
        \delta = \frac{\theta^{+}-\theta^{-}}{\sqrt{\Delta U_{\max}}} \epsilon
    \end{equation}
    where $\Delta U_{\max}$ is the local maximum of the double well potential $U$. In our numerical
    solution, we choose the parameters $\epsilon$ so the interface is solved in the mesh resolution,
    $h$. Thus,
    \begin{equation*}
        \epsilon = \frac{\sqrt{\Delta U_{\max}}}{\theta^{+}-\theta^{-}} h = \frac{h}{2\sqrt{2}}.
      \end{equation*}

    \item The parameter $\gamma$ represents the energy per unit area of the Modica-Mortola
      potential, or equivalently, the cost of introducing holes in the domain. The smaller it is,
      the less energy it costs to introduce holes into the solution and, therefore, solutions with
      more and smaller holes are obtained. If $\gamma$ is set to zero, the problem is ill-posed,
      and a sequence of solutions can be found with increasingly smaller holes and higher
      potential energy. Unlike the remaining model parameters, $\gamma$ has a direct physical
      interpretation, for it is responsible for the size and number of holes, and its value
      depends on the type of solution we wish to obtain). In addition, it must be related to the
      potential energy that is maximized, since the Modica-Mortola potential must be of the order of
      the potential energy to have an effect on evolution. Taking all this into account, a reference
      $\gamma$ is proposed that employs the value of the initial potential energy and boundary measure:
    \begin{equation}
        \gamma \sim \frac{|V(\mbs{u_0};\rho_0)|}{|\partial \Omega|}
    \end{equation}

    \item The mobility $\kappa$ controls the kinetics of the phases and is related to the
      Modica-Mortola potential. Taking into account the density evolution
      equation~\eqref{eq-transport} and the definition of the chemical
      potential~\eqref{eq-chemical-pot}, it is possible to perform a dimensional analysis to estimate
      the value of this parameter. Assuming that we want to obtain the stationary value at a
      characteristic time $T$, from Eq.~\eqref{eq-transport} we have that
    \begin{equation}
        \frac{1}{T} \sim \kappa \frac{\mu_c}{L^2},
      \end{equation}
      where $\mu_c$ and $L$ are a characteristic chemical potential and length of the problem. Also
      considering Eq.~\eqref{eq-chemical-pot}, we can estimate
    \begin{equation}
        \mu_c \sim \frac{\gamma \; \overline{U'}}{\epsilon \; \overline{\hat{\rho}'}},
      \end{equation}
      where the bars are used to indicate characteristic values of the variables. 
      For example, we have selected $\overline{U'} = \frac{2\Delta
        U_{\max}}{\theta^{+}-\theta^{-}}$ and $\overline{\hat{\rho}'} = \hat{\rho}'(0) =
      k/4$. Replacing these expressions we obtain
    \begin{equation}
        \kappa \sim \frac{\epsilon \; L^2}{T  \; \gamma} \frac{\theta^{+}-\theta^{-}}{ \;2\Delta U_{\max}}\overline{\hat{\rho}'}.
      \end{equation}
      Note that $\kappa$ depends linearly on $\epsilon$, which in turn depends on the mesh size in the discrete problem.

    \item The parameter $k$ controls the transition slope between phases $0$ and $1$ in the change
      of variable through the logistic function. This parameter is purely algorithmic and thus can
      be selected by the user. Changing its value can modify the convergence of the method, but it
      should not change the final result. By analogy with the density-based Modica-Mortola
      potential, we get $k$ such that $U(\theta)$ is (almost) a translation of $U(\rho)$. This is
      obtained by ensuring that the value of the minima are $\theta^{\pm}=\pm \frac{1}{2}$ and $k$
      is such that at those points the value $\rho_{\min}$ is obtained. Hence,
    \begin{equation}
        k = 2 \log\left(\frac{\rho_{\min}-1}{\rho_{\min}}\right).
    \end{equation}
    Also, it must be taken into account that $k$ changes the magnitude of certain terms of the
    incremental potential. In particular those of the Modica-Mortola potential, since the magnitude
    of the variable $\theta$ is modified depending on the selected parameter to obtain a minimum density
    close to zero (typically $\rho_{\min} = 10^{-3}$). Hence, different values of $k$ give rise to
    different values of $\gamma$ and $\epsilon$ so that the same solutions are obtained. In
    particular, for two different $k_1$ and $k_2$, such that $k_2=c \cdot k_1$, and from Eq.~\eqref{eq-modica},
    \begin{align}
        \frac{\gamma_{k_1}}{\epsilon_{k_1}}  & = c^4  \; \frac{\gamma_{k_2}}{\epsilon_{k_2}},
          \\
        \gamma_{k_1} \cdot \epsilon_{k_1} & = c^2  \;  \gamma_{k_2} \cdot \epsilon_{k_2} .
    \end{align}

\end{itemize}

\section{Numerical examples}
\label{sec-examples}

In this section, we examine the performance of the \vto\  method, analyzing the results obtained with it in
standard benchmarks for topology optimization, and by comparing it with two other commonly employed algorithms, the SIMP
method and an Allen-Cahn phase field discretization.

The first topology optimization method employed in our comparisons is the popular SIMP approach (Solid Isotropic Material with Penalization) \cite{bendsoe2004ik, sigmund2001simp,sigmund2013ht}. In this method, the domain of the structure is discretized using finite elements, and the design variable is the relative density of each element, a variable with values between $\rho_{\min}$ and~1. In every element, the stiffness is scaled using a power-law of the relative density with an exponent that, in our computations, we take to be equal to~3. In the comparisons shown in this section, we will employ a filter that updates the element sensitivity with a weighted mean of the sensitivities of its neighboring elements. The elements of this neighborhood are those inside a ball of radius~$r_b$.

The second topology optimization method for our numerical comparison is a purely FEM approach of the Allen-Cahn phase-field
method proposed by Takezawa \emph{et al.}\cite{takezawa2010ps}, which we will refer to as
{\take}. In this method, a phase-field
function is defined on the problem domain that indicates the ratio of void and solid, taking values
between $\rho_{\min}$ and~1. Then, the optimization process is performed by integrating a time-dependent
reaction-diffusion of the phase-field function (a modified Allen-Cahn equation with homogeneous
Neumann boundary conditions) coupled to the mechanical problem. Just as the \vto, the phase field
evolves dictated by a gradient flow, but its values are favored to $\rho_{\min}$ and~1 by a double-well
potential that depends on the local potential energy sensitivity. Since Allen-Cahn's flow does not
preserve volume, this constraint is approximately enforced with a Lagrange multiplier that is iteratively updated.
This method cannot nucleate holes; therefore, the initial phase field must include enough of them. The
method requires careful selection of three parameters, namely, $\kappa$, the diffusivity constant
for the gradient flow, $\eta$, a scaling factor of the sensitivity in the double-well potential, and
$\chi$, a penalization parameter that controls the iterative update of the Lagrange multiplier and,
hence, the approximate volume preservation.

When referring to these two methods employed for the comparisons that follow, we are aware that
there are specific implementation details that have been modified in recent works,
improving the performance of their original formulations. 
Since it is impossible to compare all existing variations of these two methods, our
presentation merely attempts to highlight some of the main characteristics of their basic
implementations, as described in their original articles and implemented in our in-house simulation code.

\subsection{The MBB beam problem}
\label{sub-mbb}

The first benchmark we use for illustrating the performance of the proposed method is the MBB
beam (see, e.g., \cite{sigmund2001simp,bartels2021ch}). As described in the literature, this
topology optimization problem is posed on a rectangular domain with dimensions $6\times1$ mm$^2$ and
an elastic isotropic material with Young's modulus $E=74000$~MPa and Poisson's ratio $\nu=0.33$. The beam is bi-supported at the ends of the lower edge, and a point load of $F = 100$ N is applied at the center of its top edge. The leftmost support is pinned, while the rightmost is simply-supported. The initial density is uniform and set to $0.5$. The geometry, boundary conditions, and external load of the MBB beam are shown in Figure~\ref{fig-MBB-beam}.

\begin{figure}[t]
    \centering
    \includegraphics[width=0.8\textwidth]{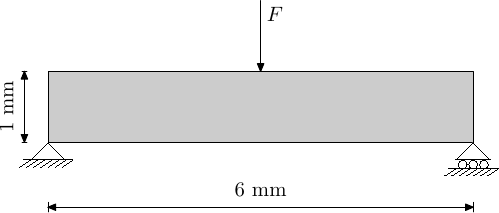}
    \caption{Geometry, supports, and loading of the MBB beam benchmark.}
    \label{fig-MBB-beam}
\end{figure}

The parameters of \vto\  are selected following the guidelines of Section~\ref{choice-parameters} except for $\gamma$, which is taken to obtain results comparable to those available in the literature~\cite{bartels2021ch}. Therefore, a similar analysis to the one performed in the description of the variable change parameter is necessary. In this case, $\gamma$ is estimated considering the definition of the chemical potential~\eqref{eq-chemical-pot},
\begin{equation}
    \frac{\gamma_{\theta}}{\epsilon_{\theta}}  \sim \overline{\hat{\rho}'}  \; \frac{\gamma_{\rho}}{\epsilon_{\rho}}, 
    \qquad
     \gamma_{\theta} \cdot \epsilon_{\theta} \sim \overline{\hat{\rho}'}^2  \; \gamma_{\rho} \cdot \epsilon_{\rho} ,
\end{equation}
giving $\gamma_{\theta} = \overline{\hat{\rho}'} \gamma_{\rho} \approx 5.77$ N/mm, when $\gamma_{\rho} =
0.9$ N/mm in the reference. The values of the parameters employed for our solutions are collected in Table~\ref{tab:parametros-MBB}.

Figure~\ref{fig-cv-den-map-mbb} shows the density maps at different times in the solution obtained
with the \vto\ method. The initial steps are very diffusive since, by design, the double-well function
has close minima. As the two wells split apart, the phases separate from each other and a
structure similar to the one available in reference works is obtained~\cite{bartels2021ch}.

\begin{figure}[p]
    \centering
    \includegraphics[width=0.9\textwidth]{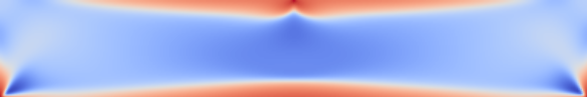}\vspace{2mm}
    \includegraphics[width=0.9\textwidth]{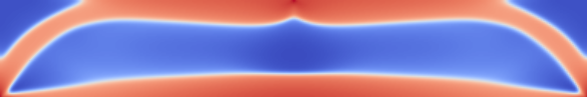}\vspace{2mm}
    \includegraphics[width=0.9\textwidth]{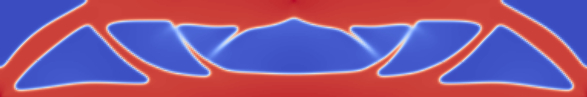}\vspace{2mm}
    \includegraphics[width=0.9\textwidth]{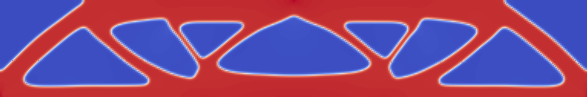}\vspace{2mm}
    \includegraphics[width=0.9\textwidth]{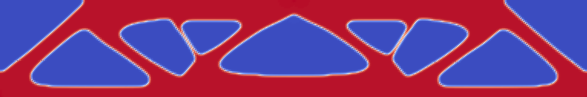}\vspace{2mm}
    \includegraphics[width=0.9\textwidth]{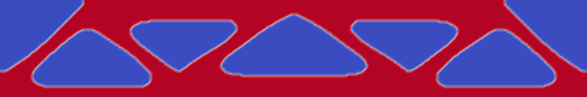}

    \vspace{5mm}
    \includegraphics[width=0.5\textwidth]{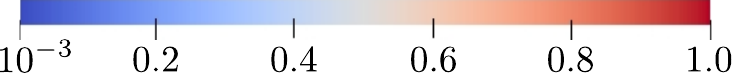}
    \caption{Relative density in the solution of the MBB problem obtained with the \vto\
      method.
      From top to bottom, snapshots correspond to $t=0.01$~s, $t=0.16$~s, $t=0.36$~s, $t=0.45$~s, $t=0.54$~s, and $t=0.7$~s.}
    \label{fig-cv-den-map-mbb}
\end{figure}

\begin{figure}[htbp!]
    \centering
    \includegraphics[width=0.8\textwidth]{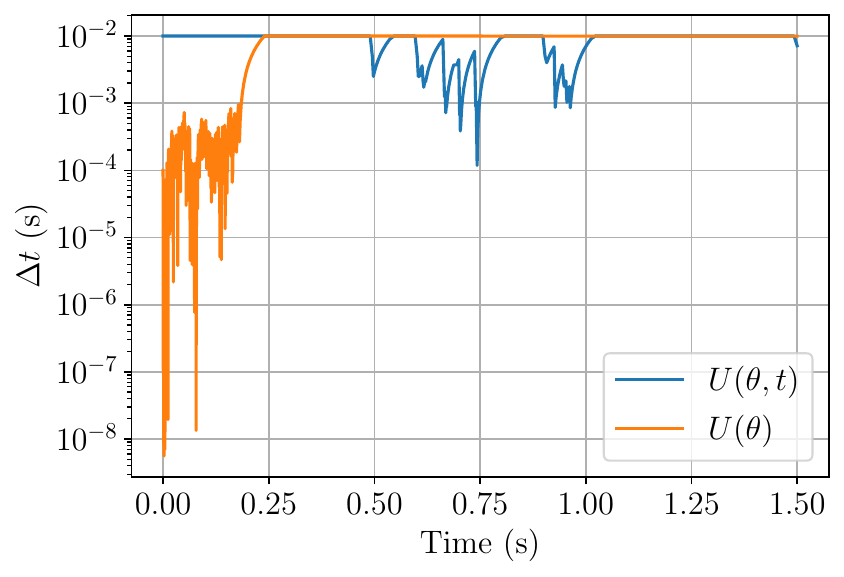}
    \caption{Time step evolution for the {\vto} method with the original double-well potential vs. the proposed continuation technique.}     \label{fig-MBB-time-step-cont}
\end{figure}

To illustrate the acceleration introduced by the continuation method, two simulations were
performed: one with the continuation method
and the other without it. In Figure~\ref{fig-MBB-time-step-cont}, the time-step size is plotted vs. time for the two simulations. For the same final simulation time ($T = 1$ s),  the simulation that
takes advantage of the continuation method can use larger time-step sizes, and hence the number of iterations and the computational time are reduced significantly.

Figure~\ref{fig-MBBenergies} depicts the evolution of the two energy terms that most contribute to
the incremental potential~\eqref{eq-incremental}, namely, the total potential energy and the
Modica-Mortola functional. The figure shows the competition between maximizing the potential energy of
the beam and minimizing the (negative) Modica-Mortola term. The continuation method employed balances these two objectives and, as shown in the figure, attains a solution that is practically stationary at $t\approx T=1$ s. After this instant, the shape of the
double-well function freezes indicating that, indeed, the parameters are selected appropriately to reach the optimal
mass distribution at~$T=1$ s.
\begin{figure}[t]
    \centering
    \includegraphics[width=0.8\textwidth]{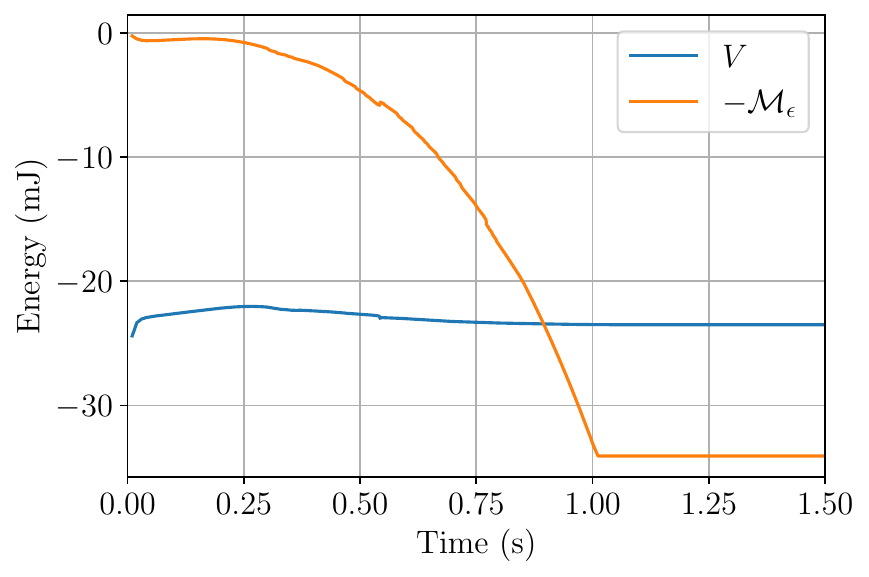}
    \caption{MBB problem. Evolution of the potential energy and Modica-Mortola functional
      in the \vto\ solution with mesh size $h = {1/64}$ mm.}
    \label{fig-MBBenergies}
\end{figure}

\begin{figure}[p]
  \centering
        \centerline{$h=1/16$ mm}
        \includegraphics[width=0.8\textwidth]{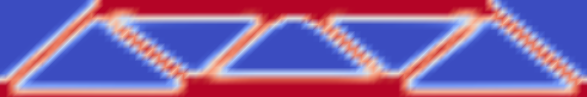}
        \includegraphics[width=0.8\textwidth]{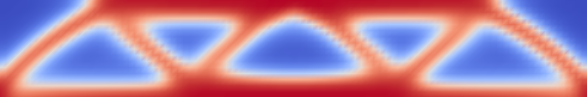}
        \includegraphics[width=0.8\textwidth]{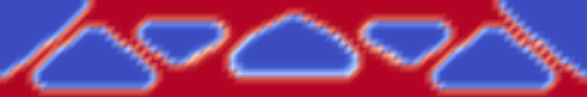}

        \vspace{5mm}
        \centerline{$h=1/32$ mm}
        \includegraphics[width=0.8\textwidth]{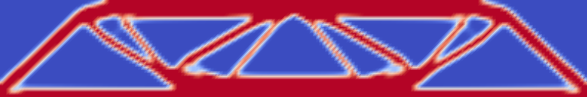}
        \includegraphics[width=0.8\textwidth]{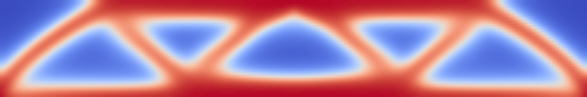}
        \includegraphics[width=0.8\textwidth]{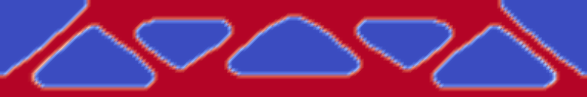}

        \vspace{5mm}
        \centerline{$h=1/64$ mm}
        \includegraphics[width=0.8\textwidth]{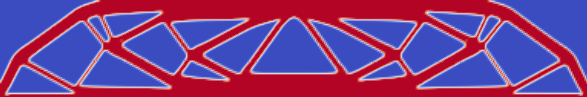}
        \includegraphics[width=0.8\textwidth]{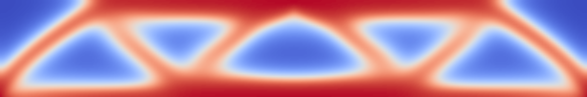}
        \includegraphics[width=0.8\textwidth]{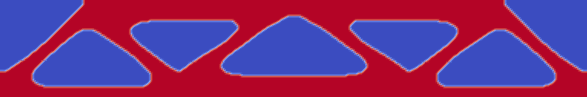}

        \vspace{3mm}
        \includegraphics[width=0.5\textwidth]{colorbar.pdf}
        \caption{Density maps of optimized structures for the MBB problem. Each block
          includes, for the indicated mesh sizes, the solution obtained with SIMP (top)
          {\take} (middle), and {\vto} (bottom) methods.}
    \label{fig-MBB-den-map-comp}
\end{figure}

\begin{figure}[t]
    \centering
    \includegraphics[width=0.8\textwidth]{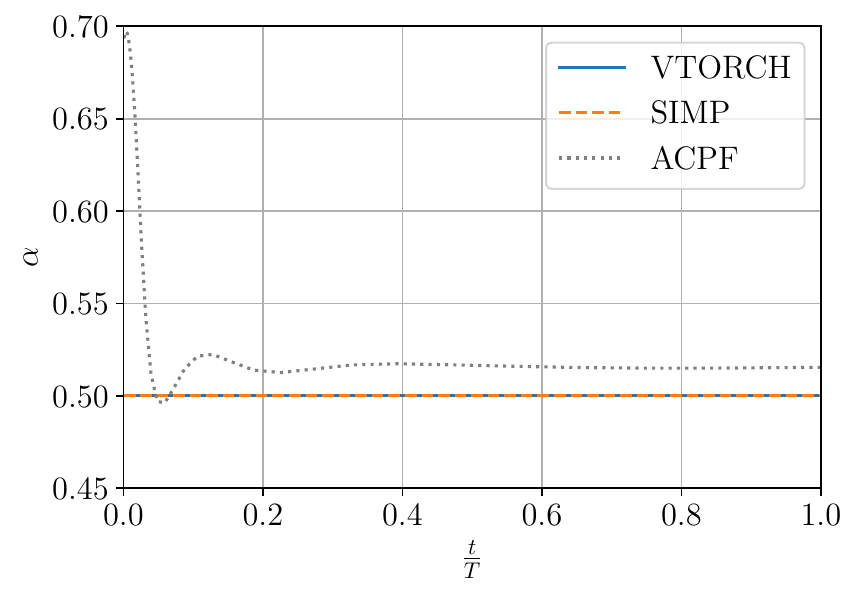}
    \caption{MBB problem. Evolution of the normalized area ($h = 1/64$~mm).}
    \label{fig-MBB-volume-comparison}
\end{figure}

Figure~\ref{fig-MBB-den-map-comp} shows the optimal topologies obtained by the
methods compared, each of them using parameters summarized in Table~\ref{tab:parametros-MBB}.
The three mass distributions qualitatively agree for the lowest resolution, although the {\take} is clearly more diffusive
than the other two. For smaller mesh sizes, {\take} and {\vto} are qualitatively similar, with {\take} still exhibiting greater diffusion. On the other hand, SIMP generates significantly more holes than the other two solutions since no length scale is imposed in the filter. It should be
recalled that the number of holes is not directly comparable between SIMP and {\vto} and ACPF
because the first method does not account for the surface energy in holes, while the other two include
it through the Modica-Mortola term. When analyzing the solutions obtained by the SIMP method,
it would seem that the area corresponding to $\rho\approx1$ is smaller than in the solutions
obtained with the other two methods. This is a postprocessing effect since mass in SIMP is
defined elementwise, and is projected to the nodes for its graphical representation, hence causing some
artificial smoothing.

Figure~\ref{fig-MBB-volume-comparison} compares the time evolution of the
normalized area, $\alpha = \frac{1}{|\Omega|}\int_{\Omega} \hat{\rho}\left(\theta(\mbs{x})\right) \;\mathrm{d}V$, for the three methods. Time is normalized with respect to the value~$T$ at which the steady-state solution is reached: in the case of the {\take} and {\vto} methods, it is always $T = 1$~s; in the case of SIMP, the evolution is considered as a function of the number of iterations. In the SIMP solution, the imposed volume is achieved up to
the desired precision by iterating. In the {\vto}, volume is exactly preserved by construction and
without the need for any special iterative scheme. Finally, in the {\take} method, the volume is not
preserved but rather remains close to the imposed value. This drift in the volume ratio is well-known
\cite{takezawa2010ps} and, although undesirable, might be compensated in practice when the solution
is post-processed with a density cut-off filter.

\begin{figure}[t]
    \centering
    \includegraphics[width=0.8\textwidth]{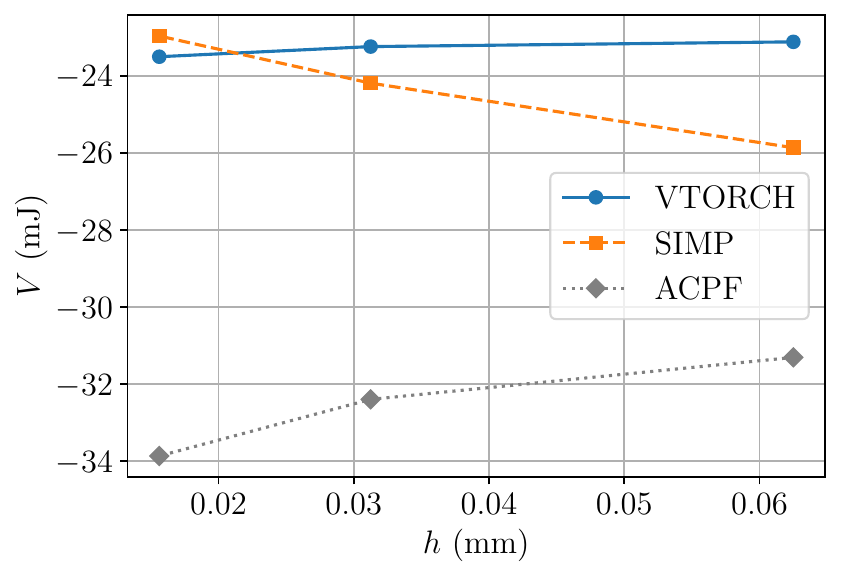}
    \caption{MBB example. Potential energy in the optimized solutions as 
     a function of the mesh size.}
    \label{fig-potential-energy-comparison}
\end{figure}

Figure~\ref{fig-potential-energy-comparison} depicts the values of the potential energy of the
optimized structures as a function of the mesh size, for the three methods compared. The results
indicate that the values are fairly independent of the mesh coarseness for the three methods and
that the {\vto} is able to obtain stiffer structures in all cases. Note, however, that the potential energy in
the SIMP and {\take} methods is defined through a penalized stiffness, which does not represent the true energy.

\begin{table}[h]
    \centering
    \begin{tabular}{lcc}
      \toprule
      SIMP \\
        $r_b/h$       & ---   & $2.0$ \\
        \midrule
      {\take} \\
        $\kappa$  &(mm$^2/$s)  & $0.1$ \\
        $\eta$    &(mm$^3$/s)       & $60$  \\
        $\chi$    &(N/mm$^5$)    & $20$ \\
        \midrule
      \vto \\
        $\epsilon/h$ &  ---    & $\frac{1}{2\sqrt{2}}$  \\
        $\gamma$  &(N/mm)      & $5.77$ \\
        $\kappa/h$  & (mm$^3$/N\,s)   & $5.07$ \\
        \bottomrule
    \end{tabular}
    \caption{Parameters used in MBB simulations.}
    \label{tab:parametros-MBB}
\end{table}

Since the functionals to be optimized in the three compared methods differ significantly (the
potential energy for {\vto} and a potential energy with penalized stiffness for SIMP and {\take})
but also in the diffusive interface and the treatment of the volume constraint, it becomes difficult
to compare the results without some post-processing. Here, we propose a uniform post-processing
approach for all methods, based on the data that would be required for manufacturing a truly
optimized structure. This step consists of analyzing new ancillary structures with density fields
that are either 1 or $\rho_0\approx 0$ everywhere, based on a given threshold. That is, for each of the methods
compared we create auxiliary solids with density fields:
\begin{equation}
\label{eq-post-processing}
\bar{\rho}_{\beta}(x) =
\begin{cases}
1 & \text{if } \rho(x) \ge \beta\ , \\
\rho_0\defined 10^{-5} & \text{otherwise}\ ,
\end{cases}
\end{equation}
and $0.5 \le \beta \le 1$. Next, for each of these models, we calculate the potential
energy under the applied loads and total mass, showing the results in Figures~\ref{fig-MBBpost-v}
and~\ref{fig-MBBpost-m}.

\begin{figure}[t]
    \centering
    \includegraphics[width=0.8\textwidth]{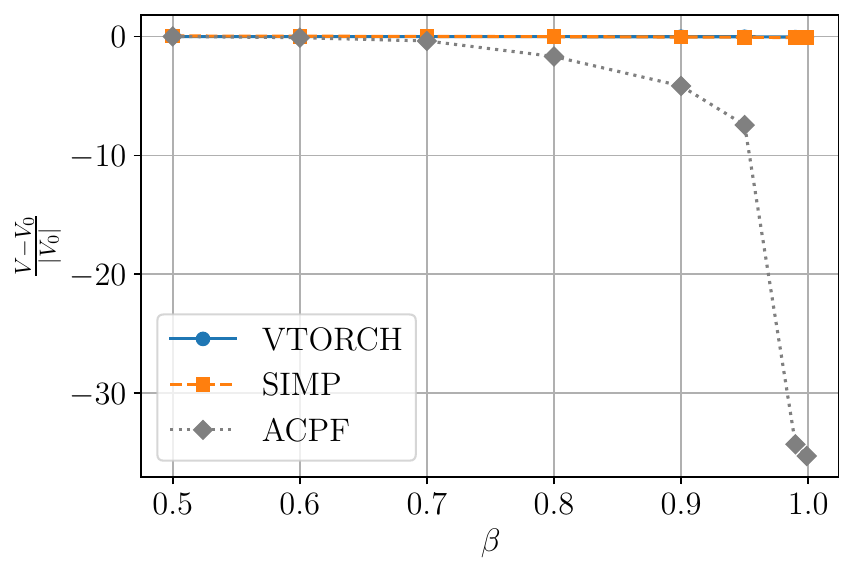}
    \caption{MBB example. Normalized potential energy in the post-processed structures
      as a function of the density threshold.}
    \label{fig-MBBpost-v}
\end{figure}

\begin{figure}[th]
    \centering
    \includegraphics[width=0.8\textwidth]{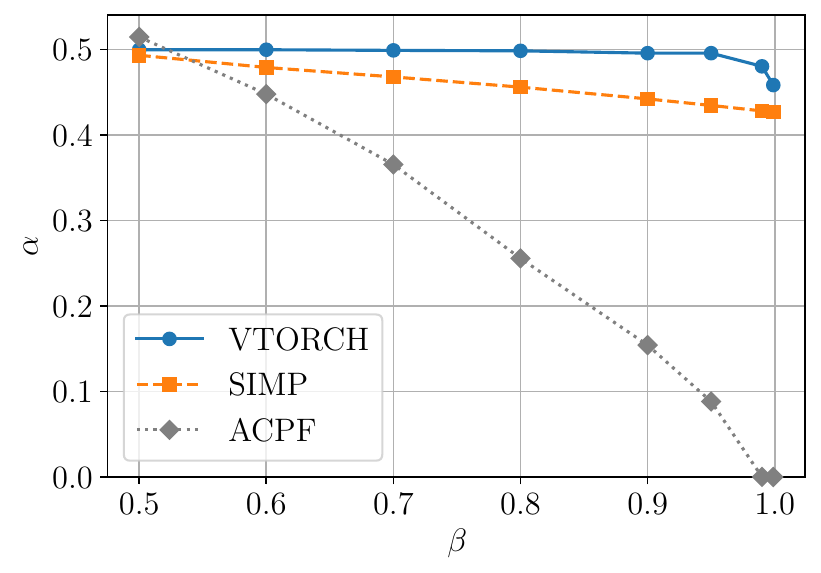}
    \caption{MBB example. Normalized area in the post-processed structures
      as a function of the density threshold.}
    \label{fig-MBBpost-m}
\end{figure}

Figure~\ref{fig-MBBpost-v} shows the results of the normalized potential
energy vs. the density threshold for the three methods. Similarly, Figure~\ref{fig-MBBpost-m}
depicts the normalized area vs. the density threshold. The potential energy is normalized
with respect to the potential energy of the {\vto} solution with $\beta=0.5$. If a method captures sharp interfaces between the
two phases, the normalized potential energy should not increase as the cut-off
threshold reaches~1 and, likewise, the normalized area should remain close to its initial value, ~$0.5$.

For the minimum threshold analyzed ($\beta=0.5$)
all three methods provide almost the same stiffness and area, with the {\take} method yielding the
maximum potential energy. This can be explained
because it is the only method of the three that does not exactly conserve the volume, and for
$\beta = 0.5$, the volume is larger than the initial one. Remarkably, even though the SIMP method obtains an
optimal shape that is qualitatively different than the optimal domains produced by
{\take} and {\vto} methods, all methods result in very similar potential energy values when $\beta = 0.5$.

Regarding the changes as the cut-off threshold $\beta$ increases, the {\vto} method is the one that
changes the least, since it produces the sharpest interfaces, showing a maximum mass reduction of
only $8\%$ when the density threshold is $\beta=0.999$. SIMP performs well and produces a sharp interface and retains volume for values $\beta \approx 1$,
showing a reduction of potential energy similar to {\vto} ($8.7\%$ and $8.4\%$, respectively).
Finally, the {\take} method provides good solutions for small values of the threshold
($\beta=0.5$). However, its high diffusivity causes its performance to deteriorate quickly, and for
values of $\beta\ge0.8$ the absolute value of the potential energy becomes more than 10 times larger than the reference one.

\FloatBarrier
\subsection{A Michell cantilever problem}
\label{sub-michell}

\begin{figure}[ht]
    \centering
    \includegraphics[width=0.6\textwidth]{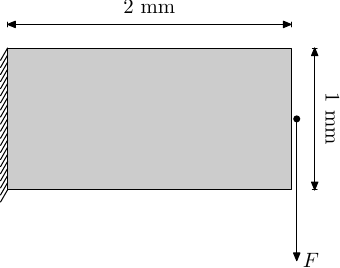}
    \caption{Geometry, supports, and loading of the Michell cantilever benchmark.}
    \label{fig-Michell_cantilever}
\end{figure}

\begin{figure}[ht]
    \centering
        \includegraphics[width=0.48\textwidth]{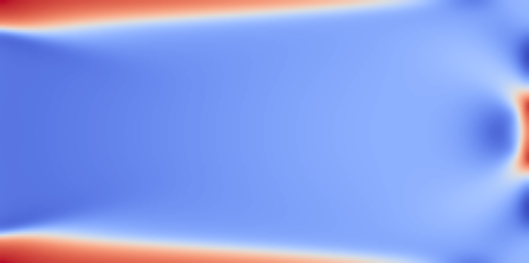}
        \includegraphics[width=0.48\textwidth]{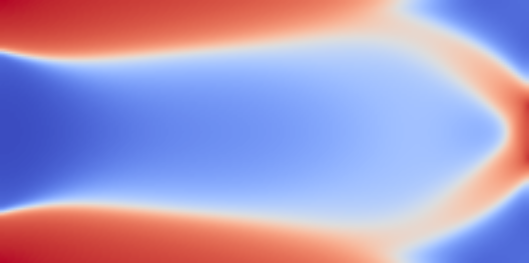}
        \includegraphics[width=0.48\textwidth]{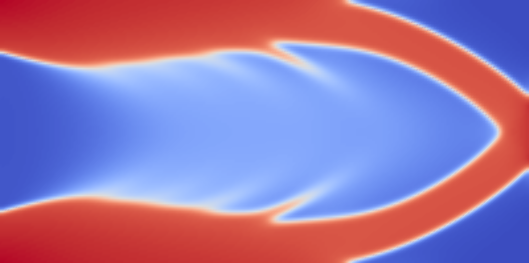}
        \includegraphics[width=0.48\textwidth]{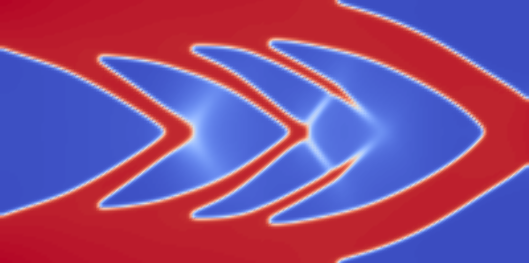}
        \includegraphics[width=0.48\textwidth]{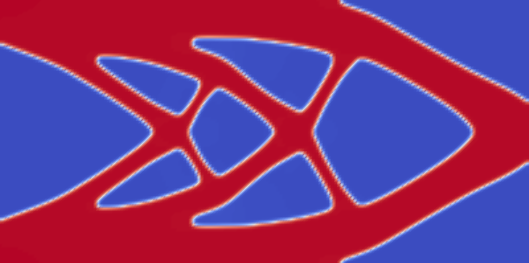}
        \includegraphics[width=0.48\textwidth]{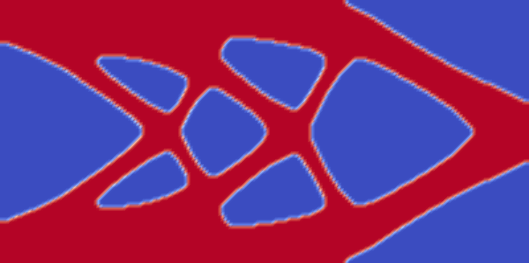}

        \vspace{3mm}
        \includegraphics[width=0.5\textwidth]{colorbar.pdf}
    \caption{Relative density maps in the solution of the Michell with the {\vto}
      method. From top to bottom, left to right, snapshots correspond to  $t=0.01$~s, $t=0.03$~s, $t=0.3$~s, $t=0.4$~s, $t=0.47$~s, and $t=0.75$~s.}
    \label{fig-cv-den-map-michell}
\end{figure}

The second selected benchmark studies a Michell cantilever (see, e.g.,
\cite{sigmund2016michell,takezawa2010ps}). The specific test we analyze is a 2D rectangular, elastic isotropic structure with dimensions $2\times1$ mm$^2$, Young's modulus $E=1$ MPa and
Poisson's ratio $\nu=0$, following~\cite{sigmund2016michell}. The cantilever is clamped at the left
edge, and a distributed load of total magnitude equal to $F=1$~N is applied
at the rightmost edge within $\pm 0.1$ mm of the symmetry axis. The initial density is uniform and
set to $0.5$. Figure~\ref{fig-Michell_cantilever} illustrates the geometry, boundary conditions, and
external load of this example. Since the parameter $\gamma$ is not defined in the references, we
have selected it following the arguments of Section~\ref{choice-parameters}, thus
\begin{equation*}
    \gamma = \frac{1}{2\,|{\partial \Omega}|} \frac{F^2 L^3}{3\rho_0EI} ,
\end{equation*}
where $L=2$ mm is the horizontal length of the rectangle, $I=\frac{1}{12}$ $\text{mm}^4$ is the cross-section inertia, $\partial \Omega = 6$ mm is the initial perimeter and $\rho_0$ the uniform density of
the initial condition. All parameters can be found in Table~\ref{tab:parametros-Michell}.

\begin{table}[h]
    \centering
    \begin{tabular}{lcc}
      \toprule
      SIMP \\
        $r_b/h$       & ---   & $2.0$ \\
        \midrule
        {\take} \\
        $\kappa$  & (mm$^2$/s)         & $5\cdot 10^{-3}$ \\
        $\eta$    & (mm$^3$/s)              & $20$ \\
        $\chi$    & (N/mm$^5$)         & $400$\\
        \midrule
        {\vto}    \\
        $\epsilon/h $   & --- & $\frac{1}{2\sqrt{2}}$ \\
        $\gamma$    & (N/mm)      & $5.33$\\
        $\kappa/h$  & (mm$^3$/Ns)  & $1.83$ \\
        \bottomrule
    \end{tabular}
    \caption{Parameters used in Michell simulations.}
    \label{tab:parametros-Michell}
\end{table}

\begin{figure}[t]
    \centering
    \includegraphics[width=0.8\textwidth]{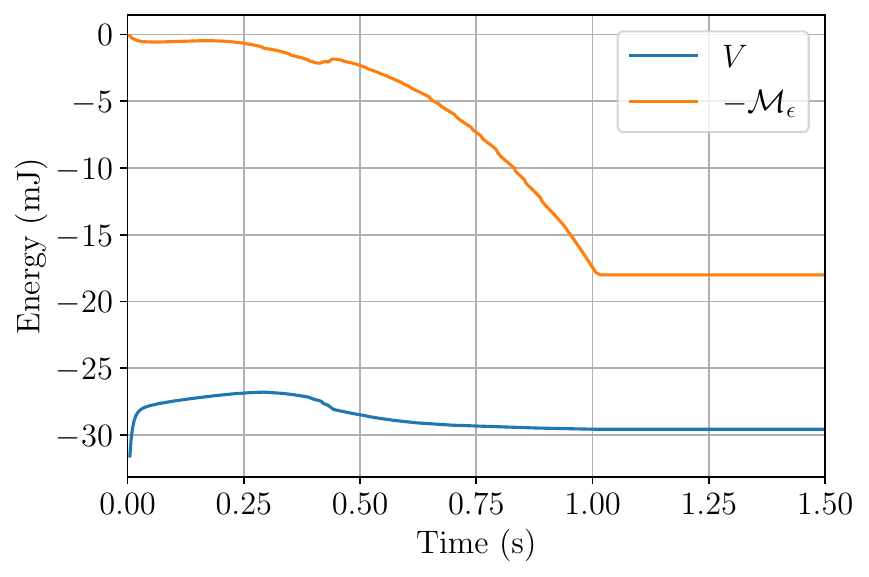}
    \caption{Michell problem. Evolution of the potential energy and Modica-Mortola functional
      in the \vto\ solution with mesh size
      ($h = 10^{-2}$ mm).}
    \label{fig-Michellenergies}
\end{figure}

Figure~\ref{fig-cv-den-map-michell} shows the density field during the optimization of the Michell
beam with the {\vto} method. As in the example of Section~\ref{sub-mbb}, the solution exhibits an
initial phase of diffuse mass separation, followed by a phase when the interfaces are sharpened,
reaching finally a stationary configuration. The value of the potential energy in the final
state is $-29.56 \; \text{N/}\text{mm}^2$, which is similar to $-29.07 \;
\text{N/}\text{mm}^2$, reported for SIMP \cite{sigmund2016michell}, although still far from the analytical value $-24.675 \; \text{N/}\text{mm}^2$ \cite{chan1962yd}. It should be noted, however,
that the last result is obtained with a perimeter constraint, and it is expected that the
{\vto} calculation will approach this value when both the mesh size and~$\gamma$ go to zero.

Figure~\ref{fig-Michellenergies} depicts
the evolution of the potential and the Modica-Mortola energies as the solution of the {\vto} method
reaches its stationary value. As in the example of the MBB beam, the competition between
minimizing $-\mathcal{M}_{\epsilon}$ and maximizing $V$ is apparent in the figure, which also
confirms that the two energies reach a saddle point, at least approximately, at time $T=1$ s.

The evolution of the normalized area can be seen in Figure~\ref{fig-Michell-volume-comparison}. In this case, $T=5$ s for {\take} method. As mentioned in the previous example, \vto\ and SIMP preserve the volume, whereas {\take} is only able to impose the constraint approximately.

\begin{figure}[t]
  \centering
  \centerline{$h=1/25$ mm}
    \includegraphics[width=0.31\textwidth]{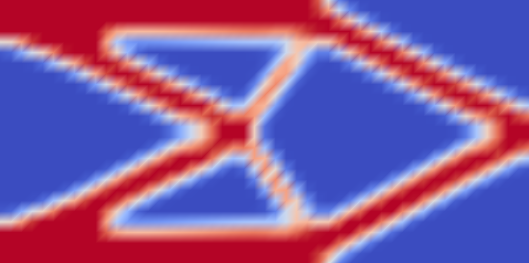}
    \includegraphics[width=0.31\textwidth]{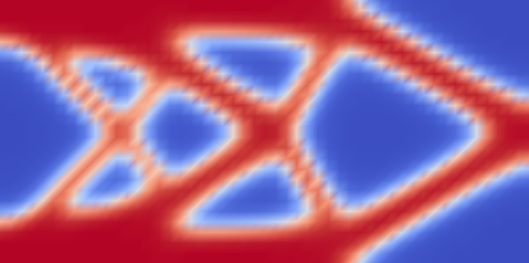}
    \includegraphics[width=0.31\textwidth]{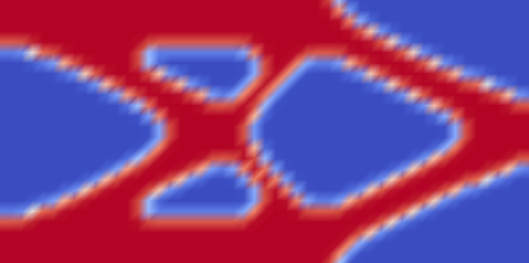}

    \centerline{$h=1/50$ mm}
    \includegraphics[width=0.31\textwidth]{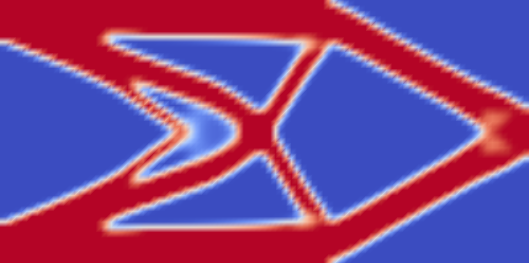}
    \includegraphics[width=0.31\textwidth]{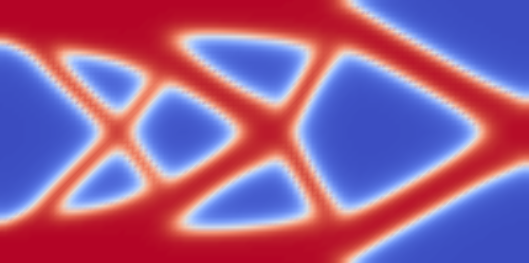}
    \includegraphics[width=0.31\textwidth]{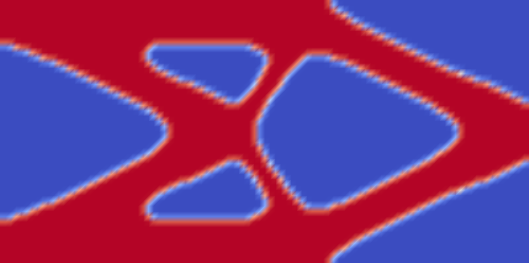}

    \centerline{$h=1/100$ mm}
    \includegraphics[width=0.31\textwidth]{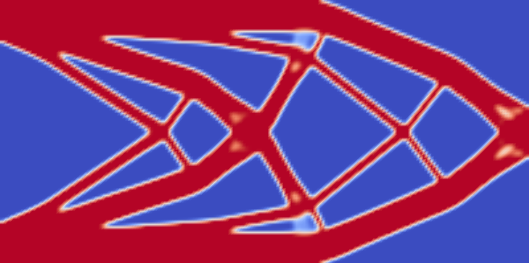}
    \includegraphics[width=0.31\textwidth]{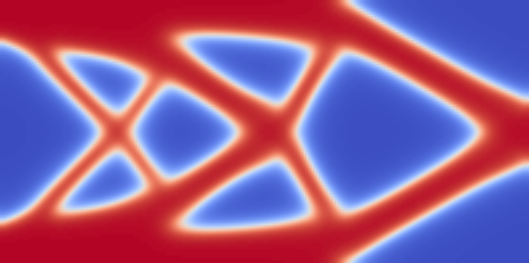}
    \includegraphics[width=0.31\textwidth]{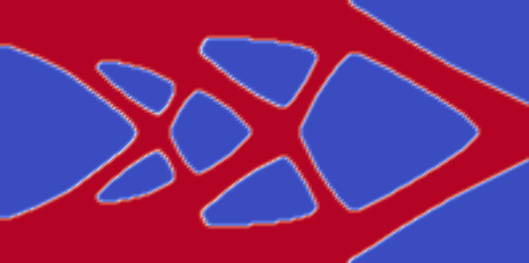}

        \vspace{5mm}
        \includegraphics[width=0.5\textwidth]{colorbar.pdf}
        \caption{Density fields of optimized structures for the Michell problem obtained
          with SIMP (left), {\take} (center), and {\vto} (right).}
    \label{fig-Michell-den-map-comp}
\end{figure}

\begin{figure}[t]
    \centering
    \includegraphics[width=0.8\textwidth]{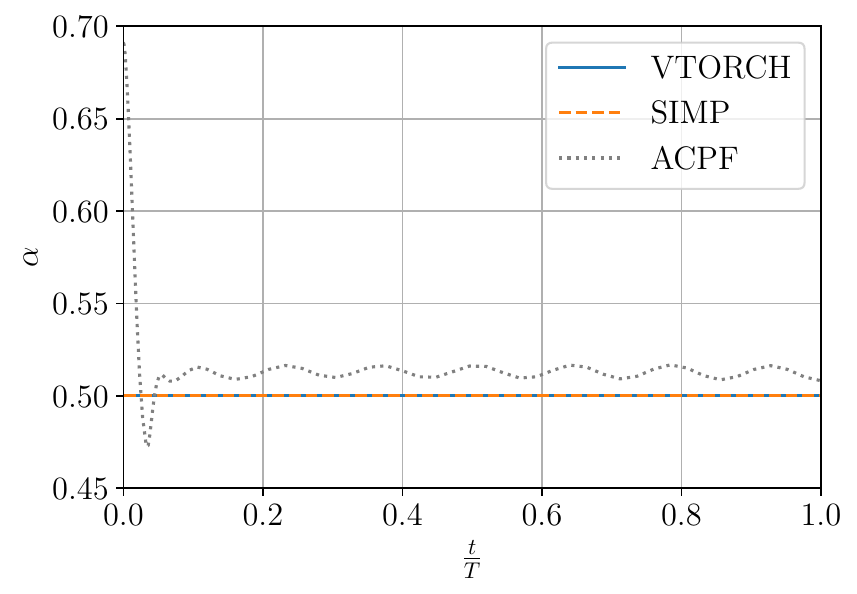}
    \caption{Michel problem. Evolution of the normalized area ($h = 10^{-2}$~mm).}
    \label{fig-Michell-volume-comparison}
\end{figure}

\begin{figure}[t]
    \centering
    \includegraphics[width=0.8\textwidth]{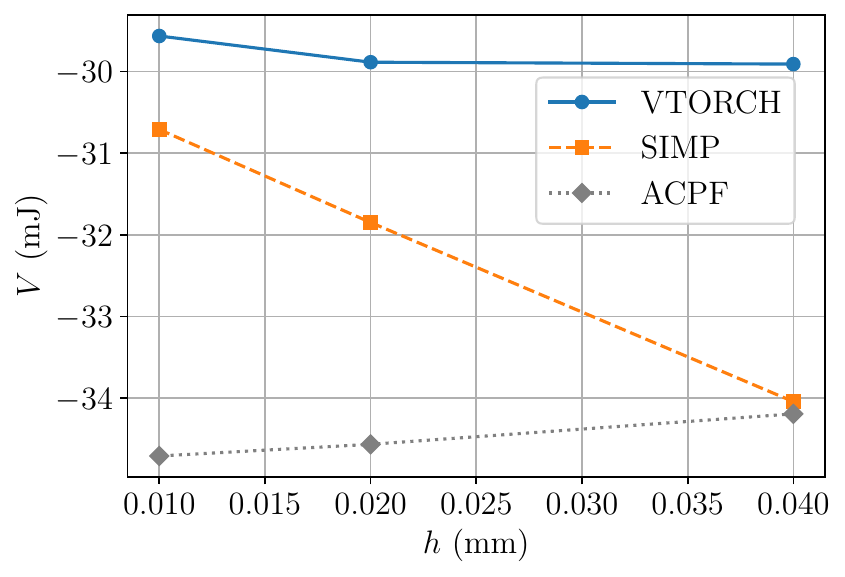}
    \caption{Michell beam. Potential energy in the optimized solutions as a function of mesh sizes.}
    \label{fig-Michell-vcomp}
\end{figure}

\begin{figure}[t]
    \centering
    \includegraphics[width=0.8\textwidth]{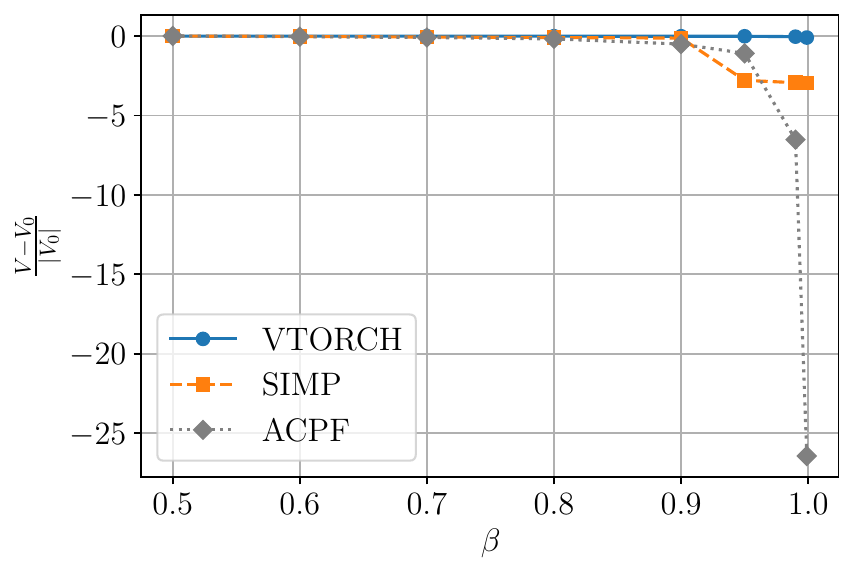}
    \caption{Michell problem. Normalized potential energy in the post-processed
      structures as a function of density threshold.}
    \label{fig-Michell-post-v}
\end{figure}

\begin{figure}[t]
    \centering
    \includegraphics[width=0.8\textwidth]{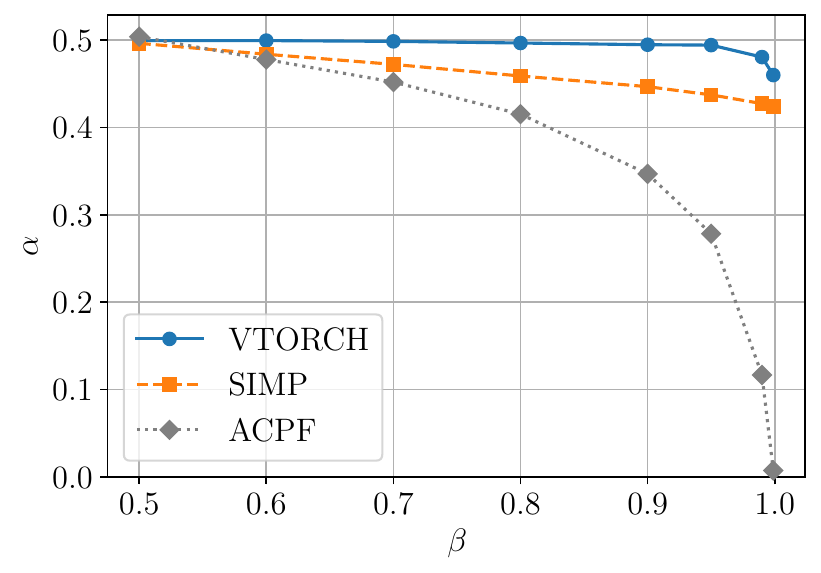}
    \caption{Michell problem. Normalized area in the post-processed
      structures as a function of density threshold.}
    \label{fig-Michell-vcomp-post-a}
\end{figure}

Figure~\ref{fig-Michell-den-map-comp} compares the optimal Michell structures as predicted by
SIMP, {\take}, and {\vto} methods. The results are qualitatively similar but, as in the MBB case, SIMP favors
a larger number of holes and connecting struts, especially for smaller mesh sizes. 

Figure~\ref{fig-Michell-vcomp} shows the potential energies of the problem
calculated in each of the formulations. In the case of {\take} and {\vto}, there are minor differences when the mesh size decreases. On the other hand, SIMP exhibits a notable increase in potential energy with mesh refinement, due to the generation of more struts. 

In Figures~\ref{fig-Michell-post-v} and~\ref{fig-Michell-vcomp-post-a}, we compare the performance of
all methods after the
post-processing explained in Section~\ref{sub-mbb}. Even though the SIMP method yields an
optimal shape that is qualitatively different than the optimal domains produced by
{\take} and {\vto} methods, very similar results in terms of potential energy are obtained for all the
methods when $\beta = 0.5$. Regarding the change with
the density cut-off $\beta$, the {\take} method behaves as in the MBB problem, exhibiting a highly
diffusive behavior and causing a rapid drop in potential energy as a consequence of the volume loss
for large values of the cut-off. SIMP loses volume linearly with $\beta$, yet it still performs well for density thresholds below $0.9$, whereas the volume loss at the interface becomes drastic when $\beta\ge0.95$, since some
of the struts are eliminated, causing a sharp drop in the potential energy, and hence the stiffness. {\vto}, produces the sharpest interface among all compared
methods and, as a result, exhibits the best behavior when $\beta\to1$. 

\FloatBarrier
\section{Summary and conclusions}
\label{sec-summary}
In this work, we have reformulated the topology optimization problem of small strain mechanics
as a Cahn-Hilliard-type evolution equation in a Hilbert space. By a suitable change of variables, we
ensure that the functional setting is not only the most natural for a finite element approximation
but also that no constraints need to be imposed on the values of the unknown density field.
More importantly, the space and time discrete problem is shown to be governed by equations that can
be obtained from the stationarity conditions of a single incremental potential, which shows that the numerical method belongs to the
class known as \emph{variational updates}. These methods are very efficient from the computational point of
view due to the guaranteed symmetry of the tangent operator.

Like other formulations in topology optimization, a double-well potential is introduced in
the variational statement of the problem to favor solid fractions that are close to~0 or~1. In our numerical solution,
we have explored the effects of a continuation method that modifies the shape of this potential as the mass
is distributed to optimize the structural stiffness. These changes have a dramatic effect on the
performance of the implementation, giving much faster convergence to the optimal topology.

We have selected some representative numerical examples that illustrate the features of the proposed method. 
First, due to the variational nature of the discrete method, only symmetric tangents are employed. Second, and due to the proposed change of
variable, no Lagrange multiplier nor special optimization strategy is required to constrain the material-to-void ratio. Third,
it is confirmed that the proposed continuation algorithm speeds up the nonlinear solver. Overall, and combined with the aforementioned continuation technique, the numerical method obtained is computationally efficient and robust.

It is difficult to compare fairly the CPU cost of different TO solution methods, and we have refrained from
specific statements in this respect. If we consider the three methods compared in Section~\ref{sec-examples}, 
they have different properties and weigh different costs in the optimization. SIMP is
the only method that does not account for surface energy and stores the density elementwise. The {\take} method, in
turn, does not solve for the pseudo-chemical potential, reducing the number of variables per node, at the expense
of losing the properties of exact mass conservation. The newly proposed {\vto} method has the largest number of
variables per node but imposes mass conservation exactly and has a symmetric tangent.

Possibly, the closest work to the current one is the published article by Bartels \emph{et al.} \cite{bartels2021ch}. We cannot make precise statements on the relative advantages of the proposed
methodology as compared to this work in terms of CPU cost, but let us note, nevertheless, that
this previous work required a model with 7 degrees of freedom per node and yielded a non-symmetric tangent. 
In turn, the {\vto} method only required 5 degrees of freedom and can use symmetric linear solvers. These two aspects amount to non-negligible CPU savings, especially for large problems. In addition, the continuation method introduced
in the current article provides for the example in Section~\ref{sub-mbb} (taken from \cite{bartels2021ch})
reduces the number of time steps of the solution from 10000, in the reference mentioned, to~400. 

Finally, let us note that the method proposed in this work can be easily adapted to topology optimization of elliptic problems besides the one described heretofore. This general class of minimization problems includes the stationary heat conduction problem, mass transport, Saint Venant warping torsion, etc. Remarkably, the isotropy of the internal energy is never exploited in the method presented so it could be extended to anisotropic problems without any essential modification.

\clearpage
\begin{acknowledgements}
  \myack%
\end{acknowledgements}

\bibliographystyle{unsrt}
\bibliography{biblio}

\end{document}